\documentclass[final]{amsart}
\usepackage{amsmath,amsfonts,amssymb,amsthm,amscd}
\usepackage[english]{babel}
\usepackage{indentfirst}
\usepackage[latin1]{inputenc}
\usepackage[T1]{fontenc}
\usepackage[mathscr]{eucal}
\usepackage{graphics}

\newtheorem{theorem}{Theorem}

\newtheorem{lemma}{Lemma}

\newtheorem{definition}{Definition}

\theoremstyle{remark}
\newtheorem{remark}{Remark}[section]

\newcommand{\beq}{\begin{equation}}
\newcommand{\eeq}{\end{equation}}
\newcommand{\ep}{epi-Lipschitzian}

\newcommand{\R}{{\mathbb R}}

\newcommand{\bd}{{\rm bd}\kern 0.12em}

\begin{document}
\title{ Representations of epi-Lipschitzian sets }
\author{Marc-Olivier Czarnecki}
\address{
Institut de Mathematiques et Modelisation de Montpellier UMR 5030 CNRS -
Universite Montpellier 2
case courier 51, Universite  Montpellier 2
Place Eugene Bataillon, 34095 Montpellier Cedex 5}
\email{marco@univ-montp2.fr}
\author{Anastasia Nikolaevna Gudovich}

\maketitle
\begin{abstract} A closed subset $M$ of a Banach space $E$ is \ep,
  i.e.,  can be represented locally as the epigraph of a Lipschitz
function, if
  and only if it is the level set of some locally Lipschitz function
  $f: E\to \R$, wich Clarke's generalized gradient does not contain
  $0$ at points in the boundary of $M$, i.e.,
such that:
\begin{eqnarray*}
 M=\{x \mid f(x)\leq 0\},\\
 0\not \in \partial f(x) \mbox{
  if } x\in \bd  M.
\end{eqnarray*}
This extends the characterization previously known in finite dimension
and answers to a standing open question
\end{abstract}
\section{\label{s1} Introduction and results}

A subset $M$ of a Banach space $E$ is \ep~(or epi-Lipschitz) at a
point $x\in M$ if it can be represented as the epigraph of a Lipschitz
function in some neighborhood of $x$. Precisely

\begin{definition}
Let $M $ be a  subset of a Banach space $E$, and let $x\in M$. The set
$M$ is \ep~ at $x$ if and only if there exists a neighborhood $U$ of
$x$, a Banach space $F$, a Banach isomorphism $A: E\to F\times \R$ and
a map $\Phi : F\to \R$ which is Lipschitz around the $F$ component of
$A(x)$, such that
$$
M\cap U=U\cap A^{-1} (epi \Phi),
$$
where $ epi \Phi=\{(\xi, t)|\Phi(\xi)\leq t\}.$ 

The set $M$ is \ep~ if
it is \ep~ around every $x\in M$. 
\end{definition}

Rockafellar~\cite{RoCan} formally introduces  the notion of \ep~ sets
in Banach spaces, and precisely gives the above definition
in~\cite{ro}, in finite dimension. We stick to this definition which
is meant by the word ``\ep'', epigraph of Lipschitz function. The
definition can also be given at a point $x\notin M$ but then requires
local closedness of $M$ to work with.

The purpose of our paper is to give a  direct and self content proof
of the characterizion of \ep~sets as level sets of
Lipschitz functions satisfying a nondegeneracy  condition, thus
extending results previously known in finite dimension.

One of our
purposes, besides answering to a standing open question, is to provide
tools for the existence of fixed points and equilibria in the non
convex case. In finite dimension, representation as level sets
permits to use an approximation method and to obtain finer results
(\cite[Theorems 2.2 and 2.3]{CCX1}) than with direct methods. In
infinite dimension, \ep~sets already play a major role
(\cite{BK2003,CK2006}) in the fixed point theory. Of course, since non empty
\ep~sets have  non empty interiors, they are not compact in infinite
dimension and one expects extension of Schauder Theorem, and certainly
not of Fan-Browder theorem.

Before stating our main result, we recall the definition of Clarke's
generalized gradient. Let $f: E\to \R$ be locally Lipschitz, the
generalized directional derivative of $f$ at $x\in E$ in the direction
$v\in E$ (see \cite{cl1}) is defined by:
$$
f^0(x,v)=\limsup_{y\to x, t>0, t\to 0}\frac{f(y+tv)-f(y)}{t}.
$$
The generalized gradient of $f$ at $x$ is the set:
$$
\partial f(x)=\{\xi\in E^* \>| \>\forall v\in X, \> f^0(x,v)\geq
\langle \xi, v\rangle\}.
$$

\begin{theorem}\label{t1} (global representation) Let $M $ be a closed subset of a Banach space $E$.
The two following assertions are equivalent:\\
(i) $M$ is \ep,\\
(ii) there is a locally Lipschitz function $f: E\to \R$
such that:
\begin{eqnarray*}
 (rep) &\mbox{(set representation)}&M=\{x \mid f(x)\leq 0\},\\
 (nd) &\mbox{(nondegeneracy)} &0\not \in \partial f(x) \mbox{
  if } x\in \bd  M.
\end{eqnarray*}
\end{theorem}

The proof of the implication $(i)\Rightarrow(ii)$ is given by
Cwiszewski and Kryszewski \cite{cwkr02} with the signed distance
function $\Delta_M$, defined by
$$
\Delta_M=d_M-d_{E\setminus M},$$ where $d_M$
is the distance function to the set $M.$ The function $\Delta_M$ is
Lipschitz with constant one, and Theorem \ref{t1} holds considering
globally Lipschitz functions instead of locally Lipschitz functions.

\begin{theorem}\label{t2}\cite[Proposition 3.10]{cwkr02} If $M\subset E$ is \ep, then for all $x\in \bd  M,$
$0\not\in\partial \Delta_M(x).$
\end{theorem}

The implication $(ii)\Rightarrow(i)$ is given in Section~\ref{s2}.

\begin{remark}\label{r1} Theorem \ref{t1} generalizes to the infinite dimension the corresponding finite dimensional
result. The proof of the remaining implication  $(i)\Rightarrow(ii)$
is then given by
\cite[Proposition 4.1]{ccz1} which \cite[Proposition 3.10]{cwkr02} generalizes
to the infinite dimension.
\end{remark}

The classical way to prove $(ii)\Rightarrow(i)$ is as follows. Take $x\in \bd  M.$ From the continuity of $f$ we get $f(x)=0.$
Then, since $0\not \in \partial f(x),$ Clarke~\cite[Corollary 1, p. 56]{cl1}  shows that
$$N_M(x)\subset \cup_{\lambda\geq0} \lambda \partial f(x),$$
where $N_M(x)$ denotes Clarke's normal cone to $M$ at $x\in M.$ It  implies that $N_M(x)$
is pointed (i.e $N_M(x)\cap - N_M(x)=\{0\}$). But in finite dimension, Rockafellar \cite{ro} characterizes  closed \ep~sets as sets having pointed normal cone. This is no longer true in infinite dimension, and Rockafellar \cite{ro} gives the following
counterexample:\\
 Let $M$ be the closed convex subset of the Hilbert space $l^2\times \R$ which is the epigraph of the function
 $$\varphi(\xi)=\sum_{j=1}^\infty j \xi_j^2, \ where\  \xi=(\xi_1, \xi_2, ...).$$
 $M$ is not \ep  \ (it has an empty interior), but 
$$
N_M(0)=\{0_{l^2}\}\times \R^-,
$$
 which is clearly pointed. Thus we cannot follow this path to prove
 the implication $(ii)\Rightarrow(i)$ in infinite dimension.

\begin{remark}\label{r2} Theorem \ref{t1} gives an answer  to
 questions raised by Bader and Kryszewski \cite{BK2003},
and by Cwiszewski and Kryszewski \cite{CK2006}. In  \cite{BK2003}, the authors raise  the question of equivalence of the then following assertions, taking $M$
to be a closed subset of a Banach space $E:$
\begin{eqnarray*}
(i) &&M \mbox{ is \ep};\\
(ii) && \forall x\in \bd  M,\qquad N_M(x)  \mbox{ is pointed}; \\
(iii) &&\forall x\in \bd  M,\qquad 0\not\in\partial \Delta_M(x) \mbox{ and }T_M(x)=\partial {\Delta_M(x)}^-,
\end{eqnarray*}
where $T_M(x)$ denotes Clarke's  cone to $M$ at $x\in M$ and $\partial {\Delta_M(x)}^-$ is the polar cone to  the
(Clarke) generalized gradient of $\Delta_M$ at $x.$

They notice the implications $(i)\Rightarrow(ii)$  and
$(iii)\Rightarrow(ii)$. They prove  the implication
$(i)\Rightarrow(iii)$ and wonder if  the implications $(ii)\Rightarrow(i),$
$(ii)\Rightarrow(iii)$ hold. The answer is no, in view of the
counterexample of Rockafellar. But  the implication $(iii)\Rightarrow(i)$
holds, from Theorem \ref{t1}.
In \cite[Example 3.2]{CK2006}, the authors take Assumption $(ii)$ to
be the definition of an \ep~set.
Theorem \ref{t1} shows that  this is consistent with the classical definition.
\end{remark}

\begin{remark} It is sometimes necessary to consider sets which are
  only \ep~on some part.  Theorem~\ref{t1} can be thus formulated more
  generally as follows.
\end{remark}

\begin{theorem}\label{t3} (local representation) Let $M $ be a subset
  of a Banach space $E$, and let $M_L\subset M$. 
The two following assertions are equivalent:\\
(i) $M$ is \ep~on $M_L$, i.e., at every point of $M_L$;\\
(ii) there is an open subset $U$ containing $M_L$ and a locally Lipschitz function $f: U\to \R$
such that:
\begin{eqnarray*}
 (rep) &\mbox{(set representation)}&M\cap U=\{x\in U \mid f(x)\leq 0\},\\
 (nd) &\mbox{(nondegeneracy)} &0\not \in \partial f(x) \mbox{
  if } x\in \bd  M\cap U.
\end{eqnarray*}
\end{theorem}

We let the reader check that the proof of Theorem~\ref{t1} (the proof
given below and the proof  of \cite[Proposition 3.10]{cwkr02}) allows to
prove Theorem~\ref{t3}, and that, with an Urysohn type argument,
Theorem~\ref{t3} implies Theorem~\ref{t1}.

\section{\label{s2} Proof of Theorem~\ref{t1}}

In this section, we prove the implication $(ii)\Rightarrow(i).$ Our proof is a nontrivial adaptation of the finite
dimensional  proof of \cite[Proposition 4.4]{ccz1}. In this paper, the proof was "left to the reader." The same result appears in \cite{Cz_these}
as Proposition 5.1, page 52, but the proof makes unnecessary use of the finite dimensional structure, among which the
scalar product.

For an element $v\in E$ and a subset $U$ of $E,$ we define (as in \cite{Cz_these} and \cite{ccz1}) the function $\lambda_{v,U}:E
\to \R\cup \{-\infty, \infty\}$ by
$$\lambda_{v,U}(x)=\inf\{t \mid x+tv\in M\cap U\}.$$
We assume $(ii)$ and consider $x\in \bd  M.$ We now find sets $U$ and $V,$ neighborhoods of $x,$ such that the function
$\lambda_{v,U}$ is real valued and Lipschitzian on $V,$ which allows us to write $M\cap V$ as the epigraph of $\lambda_{v,U}$
restricted to a certain subspace.

Since
$0\not\in\partial f (x),$ by the definition of $\partial f (x),$ and from the positive homogeneity of the generalized directional derivative, 
there exists $v\in E$ such that $\|v\|=1$ and
$$f^0(x,v)<0.$$

Take a real number $\alpha>0$ such that
$$f^0(x,v)<-\alpha.$$
By definition of the generalized directional derivative $f^0(x,v),$ there exists a real number $r>0$ such that
\begin{equation}\label{e1}
 \forall y\in B(x,2r),\qquad \forall t\in (-2r,2r)\setminus\{0\},\qquad
\frac{f(y+tv)-f(y)}{t}<-\alpha .
\end{equation}
Strictly speaking, the definition of  $f^0(x,v)$ yields the above
inequality for positive real numbers $t$, but it is easily deduced for
negative $t$.\footnote{ Assume that (\ref{e1}) to be valid for positive
  $t$, with constant $4r$. If $-2r<t<0$, then  $ y+tv \in B(x, 4r)$,
  $(f(y+tv-tv)-f(y+tv))/(-t)<- \alpha$, i e., $(f(y+tv)-f(y))/t<-\alpha$.}
Without loss of generality, possibly taking a smaller  real number
$r>0$, 
we can assume that the function $f$ satisfies  the Lipschitz condition  on $B(x,r)$ with a
some constant $k.$ Set 
$$
\varepsilon=inf\left\{\frac{r}{4},\frac{\alpha r }{4k}\right\}.
$$
\begin{lemma}\label{l1} The function $\lambda_{v,B(x,r)}$ is real-valued on $B(x,\varepsilon)$ and
\beq \label{e8}
\forall y\in B(x,\varepsilon),\qquad
\mid\lambda_{v,B(x,r)}(y)\mid \leq\frac{r}{4}.\eeq
\end{lemma}

{\bf Proof of Lemma~\ref{l1}.} Take $y\in B(x,\varepsilon),$ the set 
$$
\{t\mid y+tv \in M\cap B(x,r)\}
$$ is bounded from below, since
the set $M\cap B(x,r)$ is bounded (write $|tv|\leq
|y+tv-x|+|x-y|<r+\varepsilon$) and it contains the real number
$\frac{r}{4},$ thus showing  that
$\lambda_{v,B(x,r)}(y)\in \R.$

Indeed, first note
$$y+\frac{r}{4}v\in B(x,r)$$
from the inequality $\|y+\frac{r}{4}v -x\|\leq \|y-x\|+\frac{r}{4}<\varepsilon+\frac{r}{4}<r.$
Using  (\ref{e1}) we get
$$f\left(y+\frac{r}{4}v\right)-f(y)<-\alpha\frac{r}{4}.$$
But 
$$
f(y)-f(x)\leq k\|x-y\|
$$ 
and 
$$f(x)=0 \mbox{ (since } x\in \bd  M).
$$
Thus
 $$f\left(y+\frac{r}{4}v\right)<-\alpha\frac{r}{4} +k\varepsilon \leq 0.$$
Hence
$$y+\frac{r}{4}v\in M$$
and 
$$\frac{r}{4}\in \{t\mid y+tv\in M\cap B(x,r)\}.$$
Then $\lambda_{v,B(x,r)}(y)\in \R$ and
$$\lambda_{v,B(x,r)}(y)\leq \frac{r}{4}.$$
Now take a real number $t<0$ such that  $y+tv\in M\cap B(x,r).$ Since $\mid t \mid = \|tv\|\leq
\|y+tv-x\|+\|y-x\|<r+\varepsilon,$ we have $\mid t \mid <2r.$ From (\ref{e1})
$$f(y+tv)-f(y)>-\alpha t.$$
But
$$f(x)-f(y)\leq k\|x-y\|.$$
Recalling that $f(x)=0$ (since $\bd  M \subset \{x\mid f(x)=0 \}$ by continuity of $f$) and that $f(y+tv)\leq 0 $
(since $y+tv \in M$), we deduce $-k\varepsilon -\alpha t \leq 0.$ Thus, $t\geq - \frac{r}{4}$ and
$$\lambda_{v,B(x,r)}(y)\geq - \frac{r}{4}.$$
\qed

From Hahn Banach theorem, there exists a continuous  linear map $\varphi$ such that $\varphi(v)=1$ and $\|\varphi\|=1.$
Set
 $$F=Ker \ \varphi.$$
 Then 
$$E=F\oplus \R v$$
 and $F$ is closed. Define  the (continuous) projection on $F$ along
 $v$ 
\begin{eqnarray*}
\pi\>:\> E&\to& F\\  
y &\mapsto &y - \varphi(y)v
\end{eqnarray*}

\begin{lemma}\label{l2} For every $y\in B(x,\varepsilon)$ and $\lambda\in \R$
\begin{equation}\label{e2}
\lambda_{v,B(x,r)}(y+\lambda v )=\lambda_{v,B(x,r)}(y)-\lambda.
\end{equation}
Moreover, the function $\lambda_{v,B(x,r)}$ is real valued on the cylinder
$$B_F(\pi(x), \varepsilon)+\R v,$$
where $B_F(a, \varepsilon )=\{b\in F \mid \|b-a\|<\varepsilon \}$,
and for every $y$ in the cylinder,
$$
y+(\varphi(x)-\varphi(y))v\in B(x,\varepsilon).
$$
\end{lemma}
{\bf Proof of Lemma~\ref{l2}.} Using definition of $\lambda_{v,B(x,r)}$ we get
\begin{eqnarray*}
\lambda_{v,B(x,r)}(y+\lambda v )&=&inf\{t\mid y+\lambda v +tv \in M\cap B(x,r)\}\\
&=&inf\{t\mid y+tv \in M\cap B(x,r)\}-\lambda\\
&=&\lambda_{v,B(x,r)}(y)-\lambda.
\end{eqnarray*}

Now we must only prove that if $y\in B_F(\pi(x), \varepsilon)+\R v,$ then $y+(\varphi(x)-\varphi(y))v\in B(x,\varepsilon)$
which finishes the proof of Lemma \ref{l2} in view of (\ref{e2}). Indeed, $y=y_F+tv,$
with $y_F\in B_F(\pi(x), \varepsilon).$ By uniqueness of
decomposition, we have
\begin{eqnarray*}
y_F&=&\pi(y)=y-\varphi(y)v\\ 
t&=&\varphi(y)
\end{eqnarray*}
and $\|y-\varphi(y)v+\varphi(x)v - x\|=\|y_F-\pi(x)\|<\varepsilon.$ \qed

\begin{lemma}\label{l3} For every $y\in B_F(\pi(x), \varepsilon)+\R v$ $$y+\lambda_{v,B(x,r)}(y)v\in \bd  M\cap B\left(x, \frac{r}{2}\right).$$
\end{lemma}

{\bf Proof of Lemma~\ref{l3}.} From the definition of $\lambda_{v,B(x,r)}$ it follows that $\lambda_{v,B(x,r)}(y)=\lim_{n\to \infty}t_n $
with $y+t_nv\in M\cap B(x,r),$ and since the set $M$ is closed, we have $$y+\lambda_{v,B(x,r)}(y)v\in M.$$

\noindent Now we show that 
$$
y+\lambda_{v,B(x,r)}(y)v\in B\left(x, \frac{r}{2}\right).
$$ Indeed,
\begin{eqnarray}\label{e3}   \| y+ \lambda_{v,B(x,r)}(y)v - x
  \|&=&\|y-\varphi(y)v-(x-\varphi(x) v)+(\lambda_{v,B(x,r)}(y)+\varphi(y)-\varphi(x))\|\nonumber\\
&=&\|\pi(y)-\pi(x)+\lambda_{v,B(x,r)}
\big(y+(\varphi(x)-\varphi(y))v\big)v\|\nonumber\\&  < &\varepsilon+\mid \lambda_{v,B(x,r)}
\big(y+(\varphi(x)-\varphi(y))v\big)\mid.
\end{eqnarray}
From Lemma \ref{l2} it follows that  $y+(\varphi(x)-\varphi(y))v\in B(x, \varepsilon),$ so Lemma \ref{l1} implies that
$$
\mid \lambda_{v,B(x,r)}
(y+(\varphi(x)-\varphi(y))v)\mid\leq \frac{r}{4}.
$$
Hence,
$$\| y+ \lambda_{v,B(x,r)}(y)v - x <\varepsilon)+ \frac{r}{4}\leq  \frac{r}{2}.$$
The proof of Lemma \ref{l3} is finished by noticing (the converse
easily leads to a contradiction) that
$$y+ \lambda_{v,B(x,r)}(y)v\not\in int M. \qed$$

\begin{lemma}\label{l4} The function $\lambda_{v,B(x,r)}$ is Lipschitz on the cylinder
$B_F(\pi(x),\varepsilon)+\R v.$
\end{lemma}
{\bf Proof of Lemma \ref{l4}.} Take $y$ and $z$ in the cylinder. Set
$$y'=y+(\varphi(x)-\varphi(y))v,$$
$$z'=z+(\varphi(x)-\varphi(z))v.$$ 
From Lemma \ref{l2}, the elements $y'$ and $z'$ belongs to $B(x,\varepsilon).$
 Hence  
$$
\|z'-y'\|\leq 2 \varepsilon,
$$ and, in view of Lemma \ref{l1},
$$
\lambda_{v,B(x,r)}(y')\leq \frac{r}{4}.
$$
 Consequently,\begin{eqnarray*}
 y'+\left(\lambda_{v,B(x,r)}(y')+\frac{k}{\alpha}\|z'-y'\|\right)v&\in& B(x,r),\\
 z'+\left(\lambda_{v,B(x,r)}(y')+\frac{k}{\alpha}\|z'-y'\|\right)v&\in&
 B(x,r).
\end{eqnarray*}
 Further, since $f$ is $k$-Lipschitz on $B(x,r)$, we get
 \begin{multline}\label{e4}f\left(z'+\left(\lambda_{v,B(x,r)}(y')+\frac{k}{\alpha}\|z'-y'\|\right)v\right)\\
- f\left(y'+\left(\lambda_{v,B(x,r)}(y')+\frac{k}{\alpha}\|z'-y'\|\right)v\right)\leq k\|z'-y'\|.
 \end{multline}
The point $y'$ belongs to the cylinder, and by Lemma~\ref{l3},
 $$
y'+\lambda_{v,B(x,r)}(y') \in  \bd  M\cap B(x,r).
$$
On the other
 hand, $\frac{k}{\alpha}\|z'-y'\|\leq r,$ and by (\ref{e1}) we obtain
  \begin{multline}\label{e5}f\left(y'+\left(\lambda_{v,B(x,r)}(y')+\frac{k}{\alpha}\|z'-y'\|\right)v\right)\\
 \leq f\left(y'+\lambda_{v,B(x,r)}(y')\right)- k\|z'-y'\|\leq -k\|z'-y'\|.
 \end{multline}
 Adding (\ref{e4}) and (\ref{e5}) we get
 $$f+\left(z'+\left(\lambda_{v,B(x,r)}(y')+\frac{k}{\alpha}\|z'-y'\|\right)v\right)\leq 0.$$
 Thus, 
$$z'+\left(\lambda_{v,B(x,r)}(y')+\frac{k}{\alpha}\|z'-y'\|\right)v\in M.$$
 Since
 $z'+\left(\lambda_{v,B(x,r)}(y')+\frac{k}{\alpha}\|z'-y'\|\right)v\in
 M\cap B(x,r),$ from the definition of $\lambda_{v,B(x,r)}$ we get
 $$\lambda_{v,B(x,r)}(z')\leq \lambda_{v,B(x,r)}(y')+\frac{k}{\alpha}\|z'-y'\|.$$
 But, in view of Lemma \ref{l2},
 \beq \label{e6}\lambda_{v,B(x,r)}(y')=\lambda_{v,B(x,r)}(y) + \varphi(y)-\varphi(x),
 \eeq
 \beq \label{e7}\lambda_{v,B(x,r)}(z')=\lambda_{v,B(x,r)}(z) + \varphi(z)-\varphi(x).
 \eeq
 Note that
 $$\|z'-y'\|\leq \|z-y\|+\mid\varphi(z)-\varphi(x)\mid,$$
 and $\mid\varphi(z)-\varphi(x)\mid\leq \|z-y\|$ since $\|\varphi\|=1.$

 \noindent Now, subtracting (\ref{e6}) from (\ref{e7}) we obtain
 $$\lambda_{v,B(x,r)}(z)\leq \lambda_{v,B(x,r)}(y)+\left(1+\frac{2k}{\alpha}\right)\|z-y\|.$$
 The converse inequality is proved by exchanging $y$ and $z.$ \qed

  \begin{lemma}\label{l5} \begin{eqnarray}M\cap B(x,\varepsilon)&=&\{y\in B(x,\varepsilon)\mid
  \lambda_{v,B(x,r)}(y)\leq 0\}\nonumber\\
  &=&\{y\in B(x,\varepsilon)\mid
  \lambda_{v,B(x,r)}(y-\varphi(y)v)\leq \varphi(y)\}.\nonumber
  \end{eqnarray}
  \end{lemma}

  {\bf Proof of Lemma~\ref{l5}.} In view of Lemma \ref{l2}, 
$$\lambda_{v,B(x,r)}(y-\varphi(y)v)=\lambda_{v,B(x,r)}(y)+
  \varphi(y)$$
 and we only prove the first equality.

  \paragraph{Take} $y\in M\cap B(x,\varepsilon),$ then $$
0\in \{t\mid y+tv\in M\cap B(x,r)\}
$$
  and $\lambda_{v,B(x,r)}(y)\leq 0.$

   \paragraph{Conversely,} take $y\in B(x,\varepsilon)$ such that $\lambda_{v,B(x,r)}(y)\leq 0.$
  From Lemma \ref{l1}, $\mid~\lambda_{v,B(x,r)}(y)\mid\leq \frac{r}{4}$ and $y+\lambda_{v,B(x,r)}(y)v\in B(x,r).$
  Applying (\ref{e1}) at the point $y+\lambda_{v,B(x,r)}(y)v$ we get
  $$f(y)\leq
  f(y+\lambda_{v,B(x,r)}(y)v)-\alpha(-\lambda_{v,B(x,r)}(y))\leq 0.$$
But, by definition of $\lambda_{v,B(x,r)}(y)$ -which, in view of
  Lemma~\ref{l1} is real valued since $y\in B(x,\varepsilon)$- and
  since the set $M$ is closed, $y+\lambda_{v,B(x,r)}(y)v\in M$ and
  $f(y+\lambda_{v,B(x,r)}(y)v)=0$. Thus $f(y)\leq 0$, i.e., $y\in M.$ \qed\\

 We define the linear map 
\begin{eqnarray*}
A\>:\>E&\to& F\times \R\\ 
y&\mapsto &(y-\varphi(y)v, \varphi(y)).
\end{eqnarray*}
Clearly, A is continuous,
  invertible map and the inverse map $A^{-1} $, which is given by $A^{-1}(x,t)=x+tv$, is continuous. 
  
  \noindent Let $\Phi:F\to \R$ be a Lipschitz function such that for  all  $y\in B_F(\pi(x),\varepsilon),$
  $$ \Phi(y)=\lambda_{v,B(x,r)}(y).$$
  The proof of the implication $(ii)\Rightarrow( i)$ of Theorem \ref{t1} is finished with the following lemma.
  
  \begin{lemma}\label{l6} $M\cap B\left(x, \frac{\varepsilon}{2}\right)=B\left(x, \frac{\varepsilon}{2}\right)\cap A^{-1}(epi \Phi).$
  \end{lemma}
  {\bf Proof of Lemma~\ref{l6}.} It is a straightforward consequence of Lemma \ref{l5} and of the definition of the linear map $\Phi,$
  if we notice that, for every $y\in B(x, \frac{\varepsilon}{2})$ 
  $$\|\pi(y)-\pi(x)\|=\|y-\varphi(y)v-(x-\varphi(x)v)\|\leq 2 \|x-y\|<\varepsilon.$$
  Further, if $y\in M\cap B(x, \frac{\varepsilon}{2}) ,$  then Lemma \ref{l5} implies that 
  $\lambda_{v,B(x,r)}(\pi(y))\leq \varphi (y).$ Since $\pi(y)\in B_F(\pi(x), \varepsilon)$
  then $\Phi(y)=\lambda_{v,B(x,r)}(\pi(y)).$ Hence, 
$$
A(y)=(\pi(y),\varphi(y))\in epi \Phi,
$$
  i.e. $y\in A^{-1}(epi \Phi).$
  
  \noindent Conversely, if $y\in B(x, \frac{\varepsilon}{2})\cap A^{-1}(epi \Phi),$
  then $\pi(y)\in B_F(\pi(x), \varepsilon)$ and \\ $\lambda_{v,B(x,r)}(\pi(y))=\Phi(\pi(y))\leq \varphi(y).$
  Hence, $y\in M$ by Lemma \ref{l5}. \qed

\end{document}